%
%
%
%
%
%
%
%
%
%
%
\input amstex
\input amsppt.sty
\documentstyle{conm-p}
%
%
\message{cellular.TeX version 0.}%
\def\centertable{\leftskip=0pt plus1fill\rightskip=0pt plus1fill}

\def\begincellular#1#2\endcellular{\relax
   \begingroup
%
%
\catcode`_=11 
\ifx\forcount\undefined 
%
%
%
%
%
%
\def\forcount #1{\relax
   \def
         \for #1=##1to ##2by ##3do
            ##4%
            \endfor #1%
      {\relax
      #1=##1\relax
      \ifnum ##3>0
         \whilenot #1\ifnum ##2<#1do
            ##4%
            \advance #1 by ##3\relax
            \endwhilenot #1%
      \else
         \while #1\ifnum ##2<#1do
            ##4%
            \advance #1 by ##3\relax
            \endwhile #1%
      \fi
      }%
   \for #1%
   }%
\let\endwhilenot=\fi
\def\whilenot #1{\relax
   \def
         \whilenotloop#1 ##1do
            ##2%
            \endwhilenot #1%
      {\relax
         \expandafter\def\csname whilenotbody\string#1\endcsname{##2}%
         \expandafter\def\csname whilenotloop\string#1\endcsname
            {\relax
               ##1%
                  \let\next=\relax
               \else
                  \csname whilenotbody\string#1\endcsname
                  \expandafter\let\expandafter\next
                        \csname whilenotloop\string#1\endcsname
               \fi
               \next
               }%
         \csname whilenotloop\string#1\endcsname
         }%
   \whilenotloop#1
   }%
\let\endwhile=\fi
\def\while #1{\relax
   \def
         \whileloop#1 ##1do
            ##2%
            \endwhile #1%
      {\relax
         \expandafter\def\csname whilebody\string#1\endcsname{##2}%
         \expandafter\def\csname whileloop\string#1\endcsname
            {\relax
               ##1%
                  \csname whilebody\string#1\endcsname
                  \expandafter\let\expandafter\next
                        \csname whileloop\string#1\endcsname
               \else
                  \let\next=\relax
               \fi
               \next
               }%
         \csname whileloop\string#1\endcsname
         }%
   \whileloop#1
   }%

 \fi
\ifx\declarecount\undefined 
%
%
\def\declarecount {\allocate0\countdef}%
\def\declaredimen {\allocate1\dimendef}%
\def\declareskip  {\allocate2\skipdef}%
\def\declaremuskip{\allocate3\muskipdef}%
\def\declarebox   {\allocate4\chardef}%
\def\declaretoks  {\allocate5\toksdef}%
\def\allocate#1#2#3{\relax
   \advance\count1#1 by 1
   \ifnum\count1#1<\count20
   \else
      \errmessage{No room for \string#3!}%
   \fi
   #2#3=\count1#1
   }%

 \fi
%
%
\def\half{0.5}%
\def\by{by}%
\def\height{height}%
\def\depth{depth}%
\def\width{width}%
\def\to{to}%
\def\zeropt{0pt}%
\def\notop{\toprulewidth=\zeropt\relax}%
\def\nobottom{\bottomrulewidth=\zeropt\relax}%
\def\noleft{\leftrulewidth=\zeropt\relax}%
\def\noright{\rightrulewidth=\zeropt\relax}%
\let\x_after=\expandafter
%
%
\declaredimen\pixelwidth
\pixelwidth=1in
\divide\pixelwidth by 300                         
\declaredimen\horizontal_rule_adjust
\horizontal_rule_adjust=\zeropt
\def\sethorizontaladjustment{\horizontal_rule_adjust=\pixelwidth}%
\declaredimen\vertical_rule_adjust
\vertical_rule_adjust=\zeropt
\def\setverticaladjustment{\vertical_rule_adjust=\pixelwidth}%
%
%
\declaredimen\leftrulewidth
\declaredimen\rightrulewidth
\declaredimen\bottomrulewidth
\declaredimen\toprulewidth
%
%
\declareskip\leftborderskip
\declareskip\rightborderskip
\declareskip\bottomborderskip
\declareskip\topborderskip
\declarecount\last_column
\declaredimen\columnwidth
\declarecount\merge_columns
\declaredimen\merge_width
\declarecount\last_row
\declaredimen\rowheight
\declarecount\merge_rows
\declaredimen\merge_height
\declarecount\rowpenalty
%
%
\declaretoks\column_info
\column_info={/}%
\declaretoks\row_info
\row_info={/}%
\def\everycolumn{\leftrulewidth=0.4pt\relax
   \rightrulewidth=\leftrulewidth
   \leftborderskip=6pt plus 1fil\relax
   \rightborderskip=\leftborderskip
   \columnwidth=\zeropt\relax
   \merge_rows=0\relax
   \merge_height=\zeropt\relax
   \columnwidth=\zeropt\relax
   }%
\def\everyrow{\toprulewidth=0.4pt\relax
   \bottomrulewidth=\toprulewidth
   \topborderskip=3pt plus 1fil\relax
   \bottomborderskip=\topborderskip
   \rowheight=\zeropt\relax
   \merge_columns=0\relax
   \merge_width=\zeropt\relax
   }%
\def\get_data#1<#2{\relax
   \def\temp##1/#1/##2/##3***{\relax
      \def\temp{##2}%
      \ifnum1=0\temp
         #2={##1/#1//}%
      \else
         ##2%
      \fi
      }%
   \x_after\temp\the#2#1/1/***%
   }%
\def\add_data#1>#2#3{\relax
   \def\temp##1/#1/##2/##3***{\relax
      #2={##1/#1/##2#3/##3}%
      }%
   \x_after\temp\the#2***%
   }%
\def\add_column_number_data{\relax
   \x_after \add_data \the\column_number>\column_info
   }%
\def\get_column_number_data{\relax
   \x_after \get_data \the\column_number<\column_info
   }%
\def\add_row_number_data{\relax
   \x_after \add_data \the\row_number>\row_info
   }%
\def\get_row_number_data{\relax
   \x_after \get_data \the\row_number<\row_info
   }%
\declarebox\temp_box
\declarebox\scratch_box
\declaredimen\temp_dimen
\declaredimen\scratch_dimen
\declareskip\temp_skip
\declarecount\temp_count
\declarecount\tracingexpansions
\tracingexpansions=0
\catcode`_=8 
%
 
      #1\relax
%
%
\catcode`_=11 
%
%
\declaretoks\column_span_info
\column_span_info={}%
\declaretoks\row_span_info
\row_span_info={}%
\let\process=\relax
\declarecount\column_number
\column_number=0
%
%
\def\column#1{\relax
   \advance\column_number \by 1
   \last_column=\column_number
   \get_column_number_data
   \add_column_number_data {#1}%
   \ignorespaces
   }%
\declarecount\row_number
\row_number=0
%
%
\def\row#1{\relax
   \advance\row_number \by 1
   \message{Scanning row \the\row_number.}%
   \last_row=\row_number
   \everyrow
   \get_row_number_data
   \add_row_number_data {#1}%
   \column_number=0
   \ignorespaces
   }%
%
%
\def\blank{\relax
   \advance\column_number \by 1
   \if\column_number>\last_column
      \advance\column_number \by -1
      \column{}%
   \fi
   \everycolumn
   \get_column_number_data
   \ifnum\merge_rows>1
      \add_column_number_data {\merge_rows=0\relax}%
   \fi
   \merge_columns=0
   }%
%
%
\def\cell#1{\relax
   \advance\column_number \by 1
   \if\column_number>\last_column
      \advance\column_number \by -1
      \column{}%
   \fi
   \everycolumn
   \get_column_number_data
%
%
   \setbox\temp_box=\vbox \bgroup
      \begingroup
         \ifnum\merge_rows>0
            \advance\row_number \by -\merge_rows
            \get_row_number_data
         \fi
         \vskip \topborderskip
         \endgroup
      \hbox \bgroup
         \begingroup
            \ifnum\merge_columns>0
               \advance\column_number \by -\merge_columns
               \get_column_number_data
            \fi
            \hskip \leftborderskip
            \endgroup
         #1\vphantom{)}%
         \hskip \rightborderskip
         \egroup
      \vskip \bottomborderskip
      \egroup
%
%
   \ifnum \merge_rows>0
      \edef\temp{\process
         {\the\merge_rows}{\the\ht\temp_box}{\the\row_number}%
         \the\row_span_info
         }%
      \x_after\row_span_info\x_after=\x_after{\temp}%
      \add_column_number_data {\merge_rows=0\relax}%
   \else
%
%
      \ifdim\ht\temp_box>\rowheight
         \let\info=\relax
         \edef\temp{\the\row_number>\info
               {\rowheight=\the\ht\temp_box\relax}}%
         \let\info=\row_info
         \x_after \add_data \temp
         \rowheight=\ht\temp_box
      \fi
   \fi
%
%
   \ifnum \merge_columns>0
      \edef\temp{\process
         {\the\merge_columns}{\the\wd\temp_box}{\the\column_number}%
         \the\column_span_info
         }%
      \x_after\column_span_info\x_after=\x_after{\temp}%
      \merge_columns=0
   \else
      \ifdim\wd\temp_box>\columnwidth
         \let\info=\relax
         \edef\temp{\the\column_number>\info
               {\columnwidth=\the\wd\temp_box\relax}}%
         \let\info=\column_info
         \x_after \add_data \temp
      \fi
   \fi
   }%
%
%
\def\rcell#1{\cell{\hfill{}#1}}
\def\ccell#1{\cell{\hfill{}#1\hfill{}}}
\def\lcell#1{\cell{#1\hfill{}}}

%
%
\def\mergeright{\relax
   \advance\column_number \by 1
   \if\column_number>\last_column
      \advance\column_number \by -1
      \column{}%
   \fi
   \everycolumn
   \get_column_number_data
   \advance\merge_columns \by 1
%
%
   \ifnum\merge_rows>1
      \add_column_number_data {\merge_rows=0\relax}%
   \fi
   }%
%
%
\def\mergedown{\relax
   \advance\column_number \by 1
   \if\column_number>\last_column
      \advance\column_number \by -1
      \column{}%
   \fi
   \everycolumn
   \get_column_number_data
   \add_column_number_data {\advance\merge_rows \by 1\relax}%
   \merge_columns=0
   }%
%
%
\def\horizontalstretch#1#2#3{\relax
   \temp_count=#2\relax
   \advance\temp_count \by -#1\relax
   \edef\temp{\the\column_span_info\process{\the\temp_count}{#3}{#2}}%
   \x_after \column_span_info\x_after=\x_after{\temp}%
   \ignorespaces
   }%
\def\verticalstretch#1#2#3{\relax
   \temp_count=#2\relax
   \advance\temp_count \by -#1\relax
   \edef\temp{\the\row_span_info\process{\the\temp_count}{#3}{#2}}%
   \x_after \row_span_info\x_after=\x_after{\temp}%
   \ignorespaces
   }%
\def\noalign#1{\ignorespaces}
\catcode`_=8 
%
 
      \ignorespaces
      #2\relax                                    
%
%
\declaredimen\expansion
\edef\everycolumn{\everycolumn\expansion=\zeropt\relax}%
\edef\everyrow{\everyrow\expansion=\zeropt\relax}%
\catcode`_=11 
%
%
\def\process#1#2#3{\relax
   \last_cell=#3\relax
   \first_cell=\last_cell
   \advance \first_cell \by -#1\relax
   \span_size=#2\relax
%
%
   \gap=\span_size
   \forcount \cell_number=\first_cell to \last_cell by 1 do
      \everycell
      \get_cell_number_data
      \advance \gap \by -\cell_size
      \advance \gap \by -\expansion
      \endfor \cell_number
%
%
   \ifdim \gap>\zeropt
      \expandable_cells=#1\relax
      \advance \expandable_cells \by 1
      \trial_expansion=\zeropt
      \whilenot\search \ifdim\gap=\zeropt do
         \ifnum \expandable_cells=0
            \advance \trial_expansion \by \expansion
         \else
            \multiply \trial_expansion \by \expandable_cells
            \advance \trial_expansion \by \gap
            \divide \trial_expansion \by \expandable_cells
            \expandable_cells=0
         \fi
         \gap=\span_size
         \forcount \cell_number=\first_cell to \last_cell by 1 do
            \everycell
            \get_cell_number_data
            \advance \gap \by -\cell_size
            \ifdim \expansion>\trial_expansion
               \advance \gap \by -\expansion
            \else
               \advance \gap \by -\trial_expansion
               \advance \expandable_cells \by 1
            \fi
            \endfor \cell_number
         \temp_dimen=1sp
         \multiply \temp_dimen \by \expandable_cells
         \ifdim \gap>-\temp_dimen
            \ifdim \gap<\temp_dimen
               \gap=\zeropt
            \fi
         \fi
         \endwhilenot \search
      \forcount \cell_number=\first_cell to \last_cell by 1 do
         \everycell
         \get_cell_number_data
         \ifdim \expansion<\trial_expansion
            \let\info=\relax
            \edef\temp{\the\cell_number>\info
                  {\expansion=\the\trial_expansion\relax}}%
            \let\info=\cell_info
            \x_after \add_data \temp
            \ifnum\tracingexpansions>0
               \message{Expanded \the\cell_number}%
               \message{by \the\trial_expansion}%
               \message{from \the\cell_size}%
               \advance \cell_size \by \trial_expansion
               \message{to \the\cell_size.}%
            \fi
         \fi
         \endfor \cell_number
   \fi
   }%
\declarecount\first_cell
\declarecount\last_cell
\declaredimen\span_size
\let\expandable_cells=\temp_count
\declaredimen\trial_expansion
\let\gap=\scratch_dimen
\let\cell_number=\row_number
\let\everycell=\everyrow
\let\get_cell_number_data=\get_row_number_data
\let\cell_info=\row_info
\let\cell_size=\rowheight
\ifnum\tracingexpansions>0
   \message{Checking row expansions.}%
\fi
\the\row_span_info
\let\cell_number=\column_number
\let\everycell=\everycolumn
\let\get_cell_number_data=\get_column_number_data
\let\cell_info=\column_info
\let\cell_size=\columnwidth
\ifnum\tracingexpansions>0
   \message{Checking column expansions.}%
\fi
\the\column_span_info
\let\process=\relax
\catcode`_=8 
%
 
%
%
\catcode`_=11 
%
%
\def\column #1{\relax\ignorespaces}%
\row_number=0
\rowpenalty=0
%
%
\def\move_right_via_lastkern #1{\relax
   \temp_dimen=#1\relax
   \ifdim \lastkern>\zeropt
      \advance \temp_dimen \by \lastkern
      \unkern
   \else
   \fi
   \kern \temp_dimen
   }%
%
%
\def\row #1{\relax
   \advance \row_number \by 1
   \everyrow
   \get_row_number_data
   \advance \rowheight \by \expansion
   \ifdim \bottomrulewidth>\zeropt
      \advance \bottomrulewidth \by \horizontal_rule_adjust
   \fi
   \column_number=0
   \par
   \ifnum \rowpenalty=0
   \else
      \penalty \rowpenalty
      \rowpenalty=0
   \fi
   \noindent
   \ignorespaces
   \message{Outputting row \the\row_number.}%
   }%
%
%
\def\blank {\relax
   \advance \column_number \by 1
   \everycolumn
   \get_column_number_data
   \advance \columnwidth \by \expansion
   \advance \merge_width \by \expansion
   \move_right_via_lastkern \merge_width
%
%
   \merge_width=\zeropt
   \merge_columns=0
   \ifnum \merge_rows>0
      \add_column_number_data
            {\merge_rows=0\relax\merge_height=\zeropt\relax}%
   \fi
   }%
%
%
\def\cell #1{\relax
   \advance \column_number \by 1
   \everycolumn
   \get_column_number_data
   \advance \columnwidth \by \expansion
   \advance \merge_height \by \rowheight
   \advance \merge_width \by \columnwidth
   \ifdim \leftrulewidth>\zeropt
      \advance \leftrulewidth \by \vertical_rule_adjust
   \fi
   \ifdim \rightrulewidth>\zeropt
      \advance \rightrulewidth \by \vertical_rule_adjust
   \fi
%
%
   \begingroup
      \advance \row_number \by -\merge_rows
      \everyrow
      \get_row_number_data
      \xdef\globaltemp{\topborderskip=\the\topborderskip\relax
         \toprulewidth=\the\toprulewidth\relax
         }%
      \aftergroup \globaltemp
      \endgroup
   \ifdim \toprulewidth>\zeropt
      \advance \toprulewidth \by \horizontal_rule_adjust
   \fi
%
%
   \ifnum \merge_columns>0
      \begingroup
         \advance \column_number \by -\merge_columns
         \everycolumn
         \get_column_number_data
         \xdef\globaltemp{\leftrulewidth=\the\leftrulewidth\relax
            \leftborderskip=\the\leftborderskip\relax
            }%
         \aftergroup \globaltemp
         \endgroup
      \ifdim \leftrulewidth>\zeropt
         \advance \leftrulewidth \by \vertical_rule_adjust
      \fi
   \fi
%
%
   \setbox\temp_box=\hbox{#1}%
   \ifdim\wd\temp_box>\zeropt
      \setbox\temp_box=\hbox \bgroup
         \kern \leftborderskip
         \box\temp_box
         \egroup
      \temp_dimen=\wd\temp_box
      \advance\temp_dimen \by \rightborderskip
      \wd\temp_box=\temp_dimen
      \ifdim\wd\temp_box=\merge_width
%
%
      \else
         \setbox\temp_box=\hbox \to \merge_width \bgroup
            \hskip \leftborderskip
            #1%
            \hskip \rightborderskip
            \egroup
      \fi
%
%
      \wd\temp_box=\zeropt
      \setbox\scratch_box=\hbox{#1)}%
      \ifdim \dp\scratch_box>\dp\temp_box
         \dp\temp_box=\dp\scratch_box
      \fi
      \ifdim \ht\scratch_box>\ht\temp_box
         \ht\temp_box=\ht\scratch_box
      \fi
      \temp_dimen=\ht\temp_box
      \advance \temp_dimen \by \dp\temp_box
      \advance \temp_dimen \by \bottomborderskip
      \advance \temp_dimen \by \topborderskip
      \ifdim \temp_dimen=\merge_height
%
%
         \temp_dimen=\bottomborderskip
         \advance \temp_dimen \by \dp\temp_box
         \scratch_dimen=\rowheight
         \advance\scratch_dimen by -\temp_dimen
         \ht\temp_box=\scratch_dimen
         \raise \temp_dimen \box\temp_box
      \else 
         \setbox\temp_box=\vbox \to \rowheight \bgroup
%
%
            \advance \topborderskip \by \rowheight
            \advance \topborderskip \by -\merge_height
            \vskip \topborderskip
            \box\temp_box
            \vskip \bottomborderskip
            \egroup
         \box\temp_box
      \fi
   \fi
%
%
%
   \ifdim \toprulewidth>\zeropt
      \setbox\temp_box=\hbox \bgroup
         \temp_dimen=\merge_width
         \ifdim \half\leftrulewidth<\pixelwidth
            \kern -\pixelwidth
         \else
            \kern -\half\leftrulewidth
         \fi
         \advance \temp_dimen \by -\lastkern
         \vrule \height \half\toprulewidth
                \depth  \half\toprulewidth
                \width  \temp_dimen
         \ifdim \half\rightrulewidth<\pixelwidth
            \temp_dimen=\pixelwidth
         \else
            \temp_dimen=\half\rightrulewidth
         \fi
         \kern -\temp_dimen
         \vrule \height \half\toprulewidth
                \depth  \half\toprulewidth
                \width  2\temp_dimen
         \egroup
      \wd\temp_box=\zeropt
      \temp_dimen=\rowheight
      \advance\temp_dimen \by -\merge_height
      \ht\temp_box=\temp_dimen
      \dp\temp_box=\merge_height
      \raise \merge_height \box\temp_box
   \fi
%
%
   \ifdim \bottomrulewidth>\zeropt
      \setbox\temp_box=\hbox \bgroup
         \temp_dimen=\merge_width
         \ifdim \half\leftrulewidth<\pixelwidth
            \kern -\pixelwidth
         \else
            \kern -\half\leftrulewidth
         \fi
         \advance \temp_dimen \by -\lastkern
         \vrule \height \half\bottomrulewidth
                \depth  \half\bottomrulewidth
                \width  \temp_dimen
         \ifdim \half\rightrulewidth<\pixelwidth
            \temp_dimen=\pixelwidth
         \else
            \temp_dimen=\half\rightrulewidth
         \fi
         \kern -\temp_dimen
         \vrule \height \half\bottomrulewidth
                \depth  \half\bottomrulewidth
                \width  2\temp_dimen
         \egroup
      \wd\temp_box=\zeropt
      \dp\temp_box=\zeropt
      \ht\temp_box=\rowheight
      \box\temp_box
   \fi
%
%
   \ifdim \leftrulewidth=\zeropt
      \temp_dimen=\rightrulewidth
   \else
      \temp_dimen=\leftrulewidth
   \fi
   \ifdim \temp_dimen>\zeropt
      \setbox\temp_box=\hbox \bgroup
         \temp_dimen=\merge_height
         \advance \merge_height \by \pixelwidth
         \ifdim \leftrulewidth>\zeropt
            \kern -\half\leftrulewidth
            \vrule \height \temp_dimen
                   \depth  \pixelwidth
                   \width  \leftrulewidth
         \fi
         \ifdim \rightrulewidth>\zeropt
            \scratch_dimen=\merge_width
            \advance \scratch_dimen \by -\half\leftrulewidth
            \advance \scratch_dimen \by -\half\rightrulewidth
            \kern \scratch_dimen
            \vrule \height \temp_dimen
                   \depth  \pixelwidth
                   \width  \rightrulewidth
         \fi
         \egroup
      \wd\temp_box=\merge_width
      \ht\temp_box=\rowheight
      \dp\temp_box=\zeropt
      \box\temp_box
   \else
      \move_right_via_lastkern \merge_width
   \fi
%
%
   \merge_width=\zeropt
   \merge_columns=0
   \ifnum \merge_rows>0
      \add_column_number_data
            {\merge_rows=0\relax\merge_height=\zeropt\relax}%
   \fi
   \ignorespaces
   }%
%
%
\def\mergeright {\relax
   \advance \column_number \by 1
   \everycolumn
   \get_column_number_data
   \advance \columnwidth \by \expansion
   \advance \merge_width \by \columnwidth
   \advance \merge_columns \by 1
   \ifnum \merge_rows>0
      \add_column_number_data
            {\merge_rows=0\relax\merge_height=\zeropt\relax}%
   \fi
   }%
%
%
\def\mergedown {\relax
   \advance \column_number \by 1
   \everycolumn
   \get_column_number_data
   \advance \columnwidth \by \expansion
   \advance \merge_width \by \columnwidth
   \move_right_via_lastkern \merge_width
   \merge_width=\zeropt
   \merge_columns=0
   \advance \merge_height \by \rowheight
   \let\info=\relax
   \edef\temp{\the\column_number>\info
         {\merge_height=\the\merge_height\relax
         \advance\merge_rows \by 1\relax}}%
   \let\info=\column_info
   \x_after \add_data \temp
   \rowpenalty=10000 
   }%
\catcode`_=8 
\def\noalign#1{\relax
   \vadjust{#1}%
   \ignorespaces
   }%
%
 
      \offinterlineskip
      \parskip=\zeropt
      \ignorespaces
      #2\relax                                    
      \par
      \endgroup
   }%
%

\def\hf{\frac{1}{2}}
\def\td{\frac{1}{3}}
\def\sx{\frac{1}{16}}

\def\bQ{\bold Q}

\def\hV{\hat V}
\def\boZ{\bold Z}
\def\hboZ{\hf\bold Z}
\def\bg{\bold g}
\def\hg{\hat g}
\def\bgh{\bold{\hat g}}
\def\cA{{\Cal A}}
\def\cB{{\Cal B}}
\def\cC{{\Cal C}}

\def\cO{{\Cal O}}

\def\cV{{\Cal V}}
\def\cW{{\Cal W}}

\def\lra{\leftrightarrow}
\def\x{\times}
\def\ox{\otimes}
\def\sk1{\vskip 10pt}

\let\pr\proclaim
\let\epr\endproclaim

\topmatter
\title
Type A Fusion Rules from Elementary Group Theory
\endtitle

\rightheadtext{Type A Fusion Rules from Elementary Group Theory}

\author Alex J. Feingold, Michael D. Weiner  
\endauthor

\leftheadtext{FEINGOLD, WEINER}

\address
Dept. of Math. Sci., 
The State University of New York,
Binghamton, New York 13902-6000
\endaddress

\email alex\@math.binghamton.edu
\endemail

\address
Math. Dept.,
Pennsylvania State University, 
Altoona, Pennsylvania 16601
\endaddress

\email mdw8\@psu.edu
\endemail

\subjclass Primary 17B67, 17B65, 81T40;
Secondary 81R10, 05E10
\endsubjclass

\thanks{We wish to thank Matthias Beck, J\"urgen Fuchs and Thomas 
Zaslavsky for valuable advice.}\endthanks

\endtopmatter

\document

\noindent{\bf{1. Introduction}}
\sk1

Fusion rules play a very important role in conformal field theory \cite{Fu}, 
in the representation theory of vertex operator algebras \cite{FHL}, and in 
quite a few other areas. We are interested in understanding fusion rules from 
a combinatorial point of view, using very elementary aspects of group
theory. Several authors have also investigated fusion rules from a
combinatorial point of view \cite{BWM,BKMW,T}, but we believe our approach is 
unique. We have used this approach successfully in \cite{AFW} to obtain 
all of the $(p,q)$-minimal model fusion rules (\cite{Wa}) from elementary 
2-groups. The $(p,q)$-minimal models are a certain series of highest weight 
representations of the Virasoro algebra (\cite{KR}) which also have
the structure of a vertex operator algebra (\cite{FLM}), and modules
for it (\cite{FZ}). The Goddard-Kent-Olive coset construction \cite{GO} 
can be used to get the unitary $(p,p+1)$-minimal models from 
representations of the affine Kac-Moody Lie algebra $\bgh$ where $\bg = sl_2$ 
is the finite dimensional Lie algebra of type $A_1$. Knowing that, we were
not surprised to notice that the approach in \cite{AFW} would give the 
fusion rules for level $k$ representations of the affine $A_1$ algebra 
from elementary 2-groups. It is perhaps a bit more surprising that the 
technique worked for all of the $(p,q)$-minimal models, not just the 
unitary ones coming from affine $A_1$ coset constructions. 
Until now our approach has only given fusion rules which have 
the special property that the fusion coefficients $N_{ij}^k$ are in 
$\{0,1\}$. In order for our approach to be of greater
significance, we had to find a way to incorporate the higher
multiplicities which can occur even in the case of affine $A_2$. In this
paper we achieve that goal in such a way as to see how the simpler special
cases only have multiplicity one or zero. We cannot yet give the complete
answer for all affine type A algebras, but we present here for all levels,
our approach for ranks $N = 1$ and $N = 2$. 

The inspiration for this work on fusion rules comes from certain explicit
constructions of vertex operator algebras (VOAs), their modules, and spaces of
intertwining operators (\cite{Fe,FFR,FRW,We}) which can all be unified into
one large algebraic system. When the VOA and its modules are indexed by a
finite abelian group, $G$, and certain technical conditions are satisfied,
we call such a system a vertex operator para-algebra (VOPA). A somewhat more 
general structure, called an abelian intertwining algebra, was defined in  
\cite{DL}. In those systems 
the fusion rules are given precisely by the group algebra of $G$. This is the
situation for level one modules of an untwisted affine Kac-Moody Lie algebra 
$\bgh$, where the group $G$ is the quotient of the weight lattice by the root 
lattice of the underlying finite dimensional Lie algebra, $\bg$.  
For example, in \cite{FFR} the spinor construction of the four level one 
modules for the orthogonal affine Kac-Moody Lie algebra of type $D_4^{(1)}$ 
was used to construct a vertex operator para-algebra with $G = \boZ_2\x\boZ_2$. 
A similar phenomenon occurred in \cite{We}, where a vertex operator
para-algebra was constructed from the four level $-\hf$ modules of the
symplectic affine Kac-Moody Lie algebras of type $C_n^{(1)}$, but the group 
$G$ could be taken as either $\boZ_4$ or $\boZ_2\x \boZ_2$. 
We were searching for the appropriate generalization of this idea which would
include the VOAs and their modules coming from higher level representations
of Kac-Moody Lie algebras (\cite{FZ,TK}) and the $(p,q)$-minimal models. 
A key aspect of the generalization would have to incorporate the fusion
rules which can come from spaces of intertwining operators which are of
dimension greater than one. We first looked at 
the spinor construction of the $(3,4)$-minimal model in \cite{FRW}, 
that is, the representations of the Virasoro algebra with central element 
$c = \hf$ and with $h = 0$, $h = \hf$ and $h = \sx$, which naturally has
two copies of the $h = \sx$ module. Let us denote those modules by 
$[0]$, $[\hf]$, $[\sx]_1$ and $[\sx]_2$, respectively. 
Having two copies of $[\sx]$ allows the usual Ising model
fusion rules to be replaced by the group $\boZ_4$ or by $\boZ_2\x \boZ_2$
because there are intertwining operators labelled by independent vectors 
which act on the two $[\sx]$ modules, sending them to either $[0]$ or $[\hf]$.
This is different from the usual situation, where one intertwining operator
labelled by one vector from the one $h = \sx$ module stands for a two-dimensional
space of operators sending that $[\sx]$ module to a linear combination 
of $[0]$ and $[\hf]$. This paper and \cite{AFW} are the result of 
looking for other situations where having more copies of each irreducible
VOA module, each copy labelled by an element of a finite abelian group,
allows the fusion rules to be ``pulled apart'', showing an underlying group
algebra. 

One suspects that such an explanation of fusion rules for the higher level
Kac-Moody algebras should come from the level one para-algebra structure 
(or perhaps, abelian intertwining algebra structure). 
For example, in type $A_{N-1}$, where the weight lattice
modulo the root lattice is $\boZ_N$, the direct sum of the $N$ level one
fundamental modules should form an abelian intertwining algebra, $\cV$. 
The $k^{th}$ tensor power of $\cV$ contains all level $k$ modules, and is 
naturally an abelian intertwining algebra graded by the group $\boZ_N^k$. Then 
the fusion rules among those level $k$ modules must be strongly influenced by
that group $\boZ_N^k$, but a precise understanding of how is still lacking.
The fact that the symmetric group $S_k$ acts naturally on the $k^{th}$ tensor
power of $\cV$ is certainly relevant, and plays a crucial role in our current
paper. 

A brief summary of our results is as follows. 
Let $\hg$ be the affine algebra of type $A_{N-1}^{(1)}$ built from 
$\bg = sl_N$. Let $\lambda_1,\cdots,\lambda_{N-1}$ denote the
fundamental weights of $\bg$, and let $P^+$ 
denote the dominant integral weights of $\bg$. The irreducible 
modules $\hV_\Lambda$ for $\hg$ of level $k \geq 1$ are indexed by 
a dominant integral highest weight $\Lambda = kc + \lambda$ where 
$$\lambda = \sum_{j=1}^{N-1} a_j\lambda_j \in P^+$$ 
satisfies the ``level $k$ condition'' $\sum_{j=1}^{N-1} a_j \leq k$.
Then the irreducible modules on level $k$ are in one-to-one 
correspondence with the set of $N$-tuples of nonnegative integers
$(i_0,i_1,\cdots,i_{N-1})$ whose sum is $k$. Such an $N$-tuple
corresponds to
$$\Lambda = kc + i_1\lambda_1 + \cdots + i_{N-1} \lambda_{N-1}.$$
Fix level $k\geq 1$ and write the product in the fusion algebra 
for such modules as
$$[\lambda]\x [\mu] = \sum_{\nu\in P^+} N_{[\lambda],[\mu]}^{[\nu]}\ [\nu].$$
The distinguished module, $[0]$, is an identity element for this product, 
and for each $[\lambda]$ there is a distinguished conjugate $[\lambda^+]$
such that $N_{[\lambda],[\mu]}^{[0]} = \delta_{\mu,\lambda^+}$. 
A knowledge of the fusion coefficients, $N_{[\lambda],[\mu]}^{[\nu]}$,
is equivalent to a knowledge of the completely symmetric coefficients
$$N_{[\lambda],[\mu],[\nu]} = N_{[\lambda],[\mu]}^{[\nu^+]}.$$

Let $G = \boZ_N^k$ and let the symmetric group $S_k$ act on $G$ by 
permuting the $k$-tuples. For $a\in G$, let $[a]$ denote the 
orbit of $a$ under this action and let $\cO$ be the set of all such orbits. 
These orbits are precisely the subsets
$$P(i_0,i_1,\hdots,i_{N-1}) = \{x\in\boZ_N^k\ |\ j\hbox{ occurs exactly }
i_j\hbox{ times in }x,\ 0\leq j\leq N-1\}$$
where $(i_0,i_1,\cdots, i_{N-1})$ is any $N$-tuple of nonnegative integers 
such that $i_0+i_1+\cdots + i_{N-1} = k$. 
We now have a bijection between $\cO$ and the set of level $k$ modules. 
For $[a],[b],[c]\in \cO$ we believe the fusion coefficients
$N_{[a],[b]}^{[c]}$ have a combinatorial description in terms of the group $G$. 
The conjugate of $[c]$ is $[-c]$ and we prefer to study the symmetric coefficients 
$$N_{[a],[b],[c]} = N_{[a],[b]}^{[-c]}.$$
The combinatorial question which we were led to investigate is as follows.
For any $[a],[b],[c]\in \cO$, the group $S_k$ acts on
$$T([a],[b],[c]) = \{(x,y,z)\in [a]\x [b]\x [c] \ |\ x + y + z = 0 \}$$
which decomposes into a finite number of orbits under that action. 
Let the number of such orbits be denoted by $M([a],[b],[c])$. 
For ranks $N = 2$ and $N = 3$ we determine $M([a],[b],[c])$ and 
show how it is related to $N_{[a],[b],[c]}$. In future work we will try to give 
the relationship for all ranks. We have the following results.

\pr{Theorem} For $N = 2$, for any integral level $k\geq 1$, with 
notation as above, we have 
$$M([a],[b],[c]) = N_{[a],[b],[c]}.$$
\epr

\pr{Theorem} For $N = 3$, for any integral level $k\geq 1$, with 
notation as above, we have 
$$M([a],[b],[c]) = {N_{[a],[b],[c]} + 1 \choose 2}.$$
\epr

\newpage

\noindent{\bf{2. Background}}
\sk1

We use the definition of fusion algebra given in \cite{Fu}.
A fusion algebra $F$ is a finite dimensional commutative associative
algebra over $\bQ$ with some basis $B = \{x_a\ |\ a\in A\}$ such that the
structure constants $N_{a,b}^c$ defined by
$$x_a \cdot x_b = \sum_{c\in A} N_{a,b}^c x_c$$
are non-negative integers. It is also required that there be a distinguished
index $\Omega\in A$ with the following properties.
Define a matrix $C = [C_{a,b}]$ by $C_{a,b} = N_{a,b}^\Omega$ and 
define an associated ``conjugation'' map $\cC: F\to F$ by
$$\cC(x_a) = \sum_{b\in A} C_{a,b} x_b.$$
It is required that $\cC$ be an involutive automorphism of $F$. This implies
that $\cC^2 = I_F$ is the identity on $F$, so $C^2 = I$ is the identity
matrix. Because $0\leq N_{a,b}^c \in\boZ$, we get that either $C = I$ or
$C$ must be an order 2 permutation matrix, that is, there is a
permutation $\sigma:A\to A$ with $\sigma^2 = 1$ and
$$C_{a,b} = \delta_{a,\sigma(b)}.$$
Since $\cC$ is an automorphism, we must also
have $\cC(x_a) \cdot \cC(x_b) = \cC(x_a \cdot x_b)$, that is,
$$x_{\sigma(a)} \cdot x_{\sigma(b)} = \sum_{c\in A} N_{a,b}^c
x_{\sigma(c)}$$
which means that $N_{{\sigma(a)},{\sigma(b)}}^{\sigma(c)} = N_{a,b}^c$.
Sometimes we may write $\sigma(a) = a^+$ and call $x_{a^+}$ the conjugate of
$x_a$. One may use it to define the non-negative integers
$$N_{a,b,c} = N_{a,b}^{c^+}$$
which, by commutativity and associativity of the structure constants, are
completely symmetric in $a$, $b$ and $c$. Furthermore, one finds that
$$\eqalign{N_{\Omega,b}^c &= N_{\Omega,b,c^+} = N_{b,c^+,\Omega}\cr
&= N_{b,c^+}^{\Omega^+} = N_{b^+,c}^{\Omega}\cr
&= C_{b^+,c} = \delta_{b,c}\cr}$$
which means that $x_\Omega$ is an identity element for multiplication in
$F$,
usually written $x_\Omega = 1$. It also follows that $\Omega^+ = \Omega$.
  
In \cite{AFW} we introduced the idea of covering a fusion algebra by a 
finite abelian group and proved that the $(p,q)$-minimal model fusion
algebra, which comes from the discrete series of $0 < c < 1$ 
representations of the Virasoro algebra, can be covered by the group
$Z_2^{p+q-5}$. The basic idea, which is only set up to handle fusion
algebras whose fusion coefficients $N_{ij}^k$ are in $\{0,1\}$, is 
as follows.

\pr{Definition} Let $(G,+,0)$ be a finite abelian group and let
$G = P_0 \cup P_1 \cup ... \cup P_N$ be a partition into
$N$ disjoint subsets with $P_0 = \{0\}$.
Let $W$ be an $N$-dimensional vector space over $\bQ$ with basis
$P = \{P_0,P_1,...,P_N\}$ and define a bilinear multiplication on $W$ by the
formula
$$P_i * P_j = \sum_{k\in T(i,j)} P_k$$
where $T(i,j) = \{k\ |\ \exists a\in P_i,\exists b\in P_j,\ a+b\in P_k\}$.
We say that such a partition is associative if the product $*$ is
associative. We say that a group $G$ covers a fusion algebra if there is 
an associative partition $P$ of $G$ and a bijection
$\Phi$ between $A$ and $P$ which gives an algebra isomorphism between $F$
and $W$ such that $\Phi(\Omega) = P_0$. \epr
 
There are certainly examples of groups with partitions which are not
associative. Perhaps it is surprising that there are any interesting
examples of groups with associative partitions. As an example of a 
nonassociative partition, let $G = \boZ_5$, and let
$P_0 = \{0\}$, $P_1 = \{1,2\}$ and $P_2 = \{3,4\}$. The
multiplication table for the corresponding $W$ is as follows.
\vskip 10pt

\noindent{Table 1: A nonassociative partition for $\boZ_5$.}
\vskip 10pt

\begincellular{\centertable}
\row{}
\cell{$P_i * P_j$}\cell{$P_0$}\cell{$P_1$}\cell{$P_2$}
\row{}
\cell{$P_0$}\cell{$P_0$}\cell{$P_1$}\cell{$P_2$}
\row{}
\cell{$P_1$}\cell{$P_1$}\cell{$P_1+P_2$}\cell{$P_0+P_1+P_2$}
\row{}
\cell{$P_2$}\cell{$P_2$}\cell{$P_0+P_1+P_2$}\cell{$P_1+P_2$}
\endcellular
\vskip 10pt

We then compute
$$(P_1 * P_2) * P_2 = (P_0+P_1+P_2) * P_2 = P_2 +(P_0+P_1+P_2)+(P_1+P_2)
= P_0 + 2P_1 + 3P_2$$
and
$$P_1 * (P_2 * P_2) = P_1 * (P_1+P_2) = (P_1+P_2) + (P_0+P_1+P_2)
= P_0 + 2P_1 + 2P_2$$
so $(P_1 * P_2) * P_2 \neq P_1 * (P_2 * P_2)$.

The fusion rule tables presented here were produced by the 
computer program of Bert Schellekens, called ``Kac'', 
available from his webpage [http://norma.nikhef.nl/\break $\sim$t58/]. 
As an example of a nontrivial fusion algebra which can be covered by a
group, let the $\cW$ algebra coming from the coset construction of
$$\frac{SU(N)_r \ox SU(N)_s}{SU(N)_{r+s}}$$
be denoted by $\cW_N(r,s)$. The fusion rules for $\cW_3(1,1)$ are as follows.

\newpage

\noindent{Table 2: Fusion rules for $\cW_3(1,1)$.}
\vskip 10pt

\begincellular{\centertable}
\row{}
\cell{[a]$\x$[b]}\cell{[0]}\cell{[1]}\cell{[2]}\cell{[3]}\cell{[4]}\cell{[5]}
\row{}
\cell{[0]}
\cell{[0]}
\cell{[1]}
\cell{[2]}
\cell{[3]}
\cell{[4]}
\cell{[5]}
\row{}
\cell{[1]}
\cell{ }
\cell{[0]+[1]}
\cell{[3]}
\cell{[2]+[3]}
\cell{[5]}
\cell{[4]+[5]}
\row{}
\cell{[2]}
\cell{ }
\cell{ }
\cell{[4]}
\cell{[5]}
\cell{[0]}
\cell{[1]}
\row{}
\cell{[3]}
\cell{ }
\cell{ }
\cell{ }
\cell{[4]+[5]}
\cell{[1]}
\cell{[0]+[1]}
\row{}
\cell{[4]}
\cell{ }
\cell{ }
\cell{ }
\cell{ }
\cell{[2]}
\cell{[3]}
\row{}
\cell{[5]}
\cell{ }
\cell{ }
\cell{ }
\cell{ }
\cell{ }
\cell{[2]+[3]}
\endcellular
\vskip 10pt
     
\sk1
Note that $\{[0], [2], [4]\}$ forms a subgroup isomorphic to $\boZ_3$.
We find that $\boZ_3^2$ covers these fusion rules as follows:
$$\eqalign{
\{(0,0)\} &\lra [0] \qquad \{(1,2),(2,1)\} \lra [1] \qquad \{(1,1)\} \lra [2] \cr
\{(0,2),(2,0)\} &\lra [3] \qquad \{(2,2)\} \lra [4] \qquad \{(1,0),(0,1)\} \lra [5] 
\cr}$$

A multiplication table which illustrates exactly how $\boZ_3^2$ covers the
$W_3(1,1)$ fusion rules is as follows.
\vskip 10pt

\noindent{Table 3: $\boZ_3^2$ covering the $\cW_3(1,1)$ fusion rules.}
\vskip 10pt

\begincellular{\centertable}
\row{}
\cell{[a]$\x$[b]}\cell{(0,0)}\cell{(1,2),(2,1)}\cell{(1,1)}\cell{(0,2),(2,0)}
\cell{(2,2)}\cell{(1,0),(0,1)}
\row{}
\cell{(0,0)}\cell{(0,0)}\cell{(1,2),(2,1)}\cell{(1,1)}\cell{(0,2),(2,0)}
\cell{(2,2)}\cell{(1,0),(0,1)}
\row{}
\cell{\vbox{\hbox{(1,2)} \hbox{(2,1)}}}
\cell{ }
\cell{\vbox{\hbox{(2,1),(0,0)} \hbox{(0,0),(1,2)}}}
\cell{\vbox{\hbox{(2,0)} \hbox{(0,2)}}}
\cell{\vbox{\hbox{(1,1),(2,0)} \hbox{(0,2),(1,1)}}}
\cell{\vbox{\hbox{(0,1)} \hbox{(1,0)}}}
\cell{\vbox{\hbox{(2,2),(0,1)} \hbox{(1,0),(2,2)}}}
\row{}
\cell{(1,1)}
\cell{ }
\cell{ }
\cell{(2,2)}
\cell{(1,0),(0,1)}
\cell{(0,0)}
\cell{(2,1),(1,2)}
\row{}
\cell{\vbox{\hbox{(0,2)} \hbox{(2,0)}}}
\cell{ }
\cell{ }
\cell{ }
\cell{\vbox{\hbox{(0,1),(2,2)} \hbox{(2,2),(1,0)}}}
\cell{\vbox{\hbox{(2,1)} \hbox{(1,2)}}}
\cell{\vbox{\hbox{(1,2),(0,0)} \hbox{(0,0),(2,1)}}}
\row{}
\cell{(2,2)}
\cell{ }
\cell{ }
\cell{ }
\cell{ }
\cell{(1,1)}
\cell{(0,2),(2,0)}
\row{}
\cell{\vbox{\hbox{(1,0)} \hbox{(0,1)}}}
\cell{ }
\cell{ }
\cell{ }
\cell{ }                                    
\cell{ }
\cell{\vbox{\hbox{(2,0),(1,1)} \hbox{(1,1),(0,2)}}}
\endcellular
\vskip 10pt
             
The simplest series of fusion algebras comes from the
affine algebra $\hg$ of type $A_1^{(1)}$ built on $\bg = sl_2$.
There are $k+1$ modules $\hV_a$ for $\hg$ of level $k \geq 1$ which
can be indexed by ``spin'' $a\in \hboZ$ with $0\leq a \leq \frac{k}{2}$. The
conformal weight of the highest weight vector in $\hV_a$ is
$$\Delta_a = \frac{a(a+1)}{k+2}.$$
The $\bg$-submodule $V_a$ of $\hV_a$ generated by a highest weight vector is
of dimension $2a+1$. There is a simple formula for the decomposition of the
tensor product $V_a\otimes V_b$ into the direct sum of $\bg$-modules $V_c$,
$$V_a\otimes V_b = \sum_{|a-b|\leq c\leq a+b} V_c$$
where the sum is only taken over those $c\in \hboZ$ such that
$a+b+c\in\boZ$.
This decomposition is independent of the level $k$ and is part of the
basic representation theory of $sl_2$. It is remarkable that the fusion rules
for level $k$ are obtained simply by truncating the above summation. 
The following formula for the $A_1^{(1)}$ level $k$ fusion rules is well-known
(\cite{TK}):
$$N_{a,b}^c = 1 \quad\hbox{ if }\quad |a-b|\leq c\leq a+b, \qquad
a+b+c\in\boZ,
\qquad a + b + c \leq k,$$
and $N_{a,b}^c = 0$ otherwise. The conditions above imply that 
$c\leq \frac{k}{2}$.

We can also express these fusion rules in the following simple way.
Re-index the $k+1$
modules $\hV_a$ on level $k$ by $m = 2a+1\in \boZ$ with $1\leq m \leq k+1$.
Then $m = dim(V_a)$ and we write $\hV_a = \hV(2a+1) = \hV(m)$. Let $p =
k+2$.

\pr{Definition} For integer $p\geq 2$, we say that the triple of integers
$(m,m',m'')$ is {\bf p-admissible}
when  $0<m,m',m''<p$, the sum $m+m'+m''<2p$ is odd, and
the ``triangle'' inequalities $m<m'+m''$, $m'<m+m''$, and $m''<m+m'$ are
satisfied.\epr

Then the fusion rules are given by $N_{m,m'}^{m''} = 1$ if $(m,m',m'')$ is
p-admissible, $N_{m,m'}^{m''} = 0$ otherwise.

The proof in \cite{AFW} is easily modified to show that the 
$A_1^{(1)}$ fusion ring for level $k$ is covered by the group 
$\boZ_2^k$. 

\pr{Theorem 1} The level $k$ fusion rules for $\hg$ of type $A_1^{(1)}$ define
a fusion algebra $F$ with $A = \{m\in \boZ\ |\ 1\leq m \leq k+1\}$ with
distinguished element $\Omega$ being $m = 1$ and the conjugate of $m$ is
$m$. The the fusion algebra $F$ is covered by the elementary abelian 2-group 
$G = \boZ_2^{k}$ with partition given by
\vskip -5pt
$$P_i = \{g\in G\ |\ \hbox{exactly }i\hbox{ coordinates of }g\hbox{ are
}1\} \hbox{ for } 0\leq i\leq k.$$ \epr
\vskip -2pt 
                                                        
We illustrate how this theorem works in the following tables. We have used
a labelling scheme where on level $k$ the $k+1$ 
generators of the fusion ring are denoted by $[i]$ for $0\leq i\leq k$ and the
distinguished element $\Omega$ is $[0]$. Following each fusion rule table is a 
table showing how the group $\boZ_2^k$ covers it. 
\vskip 5pt

\noindent{Table 4: Fusion Table for $A_1$ of level $k = 2$.}
\vskip 5pt

\begincellular{\centertable}
\row{}
\cell{[i]$\x$[j]}\cell{[0]}\cell{[1]}\cell{[2]}
\row{}
\cell{[0]}\cell{[0]}\cell{[1]}\cell{[2]}
\row{}
\cell{[1]}\cell{}\cell{[0]}\cell{[2]}
\row{}
\cell{[2]}\cell{}\cell{}\cell{[0]+[1]}
\endcellular
\vskip 5pt

\noindent{Table 5: Group $\boZ_2^2$ covering the Fusion Table 
for $A_1$ of level $k = 2$.}
\vskip 5pt

\begincellular{\centertable}
\row{}
\cell{[i]$\x$[j]}\cell{(0,0)}\cell{(1,1)}\cell{(1,0),(0,1)}
\row{}
\cell{(0,0)}\cell{(0,0)}\cell{(1,1)}\cell{(1,0),(0,1)}
\row{}
\cell{(1,1)}\cell{}\cell{(0,0)}\cell{(0,1),(1,0)}
\row{}
\cell{(1,0),(0,1)}\cell{}\cell{}\cell{\vbox{\hbox{(0,0),(1,1)} \hbox{(1,1),(0,0)}}}
\endcellular

\newpage

\noindent{Table 6: Fusion Table for $A_1$ of level $k = 3$.}
\vskip 5pt

\begincellular{\centertable}
\row{}
\cell{[i]$\x$[j]}\cell{[0]}\cell{[1]}\cell{[2]}\cell{[3]}
\row{}
\cell{[0]}\cell{[0]}\cell{[1]}\cell{[2]}\cell{[3]}
\row{}
\cell{[1]}\cell{}\cell{[0]}\cell{[3]}\cell{[2]}
\row{}
\cell{[2]}\cell{}\cell{}\cell{[0]+[2]}\cell{[1]+[3]}
\row{}
\cell{[3]}\cell{}\cell{}\cell{}\cell{[0]+[2]}
\endcellular
\vskip 10pt

\noindent{Table 7: Group $\boZ_2^3$ covering the Fusion Table 
for $A_1$ of level $k = 3$.}
\vskip 10pt

\begincellular{\centertable}
\row{}
\cell{[i]$\x$[j]}\cell{(0,0,0)}\cell{(1,1,1)}\cell{(0,1,1),(1,0,1),(1,1,0)}\cell{(1,0,0),(0,1,0),(0,0,1)}
\row{}
\cell{(0,0,0)}\cell{(0,0,0)}\cell{(1,1,1)}\cell{(1,0,0),(0,1,0),(0,0,1)}\cell{(0,1,1),(1,0,1),(1,1,0)}
\row{}
\cell{(1,1,1)}\cell{}\cell{(0,0,0)}\cell{(0,1,1),(1,0,1),(1,1,0)}\cell{(1,0,0),(0,1,0),(0,0,1)}
\row{}
\cell{\vbox{\hbox{(0,1,1)} \hbox{(1,0,1)} \hbox{(1,1,0)} }}
\cell{}\cell{}\cell{\vbox{
\hbox{(0,0,0),(1,1,0),(1,0,1)} 
\hbox{(1,1,0),(0,0,0),(0,1,1)}  
\hbox{(1,0,1),(0,1,1),(0,0,0)} }}
\cell{\vbox{
\hbox{(1,1,1),(0,0,1),(0,1,0)} 
\hbox{(0,0,1),(1,1,1),(1,0,0)}  
\hbox{(1,0,0),(1,0,0),(1,1,1)} }}
\row{}
\cell{\vbox{\hbox{(1,0,0)} \hbox{(0,1,0)} \hbox{(0,0,1)} }}
\cell{}\cell{}\cell{}\cell{\vbox{
\hbox{(0,0,0),(1,1,0),(1,0,1)} 
\hbox{(1,1,0),(0,0,0),(0,1,1)}  
\hbox{(1,0,1),(0,1,1),(0,0,0)} }}
\endcellular

\vskip 15pt

The next simplest series of fusion algebras comes from the 
affine algebra $\hg$ of type $A_2^{(1)}$ built from $\bg = sl_3$. 
There are $(k+1)(k+2)/2$ irreducible modules $\hV_\Lambda$ 
for $\hg$ of level $k \geq 1$ which can be indexed by a dominant integral 
highest weight $\Lambda = kc + \lambda$ where 
$\lambda = a_1\lambda_1 + a_2\lambda_2\in P^+$ is a dominant integral weight
for the finite dimensional algebra $\bg$ satisfying the ``level $k$ condition'' 
$a_1 + a_2 \leq k$ for $0\leq a_i\in\boZ$. The irreducible $\bg$-submodule 
$V_\lambda$ of $\hV_\Lambda$ generated by a highest weight vector is
of dimension $(a_1+1)(a_2+1)(a_1+a_2+2)/2$. 
The decomposition of the tensor product of irreducible $\bg$-modules
\vskip -6pt
$$V_\lambda\otimes V_\mu = \sum_{\nu\in P^+} Mult_{\lambda,\mu}^\nu V_\nu$$
into the direct sum of irreducible $\bg$-modules includes multiplicities. 
This decomposition is independent of the level $k$ and is part of the 
basic representation theory of $sl_3$. The fusion rules for level $k$ 
are obtained by a rather subtle truncation of the above summation. The
tensor product multiplicity, $Mult_{\lambda,\mu}^\nu$, is equal to 
the multiplicity of the trivial module $V_0$ in the triple tensor product 
$V_\lambda\otimes V_\mu\otimes V_{\nu^+}$, where $V_{\nu^+}$ is the 
contragredient module of $V_\nu$. This means that $\nu^+$ is the highest
weight of the dual space of $V_\nu$, so $\nu^+$ is the negative of the
lowest weight of $V_\nu$. Then
$$Mult_{\lambda,\mu,\nu} = Mult_{\lambda,\mu}^{\nu^+}$$
is the multiplicity of $V_0$ in the triple tensor product 
$V_\lambda\otimes V_\mu\otimes V_{\nu}$. This number is clearly symmetric
in $\lambda$, $\mu$ and $\nu$. 

We have the following result from \cite{BMW} giving the $A_2$ tensor
product multiplicities and the fusion rules on level $k$. 

\pr{Theorem 2} \cite{BMW} Let 
$\lambda = a_1 \lambda_1 + a_2 \lambda_2$, 
$\mu = b_1 \lambda_1 + b_2 \lambda_2$, and
$\nu = c_1 \lambda_1 + c_2 \lambda_2$
be dominant integral weights for $\bg = sl_3$, and let
$$\cA = \td\left( 2(a_1+b_1+c_1) + a_2+b_2+c_2 \right)\ , \quad
\cB = \td\left( a_1+b_1+c_1 + 2(a_2+b_2+c_2) \right).$$
Define $k_0^{max} = \min(\cA,\cB)$ and
$$k_0^{min} = \max(a_1+a_2,b_1+b_2,c_1+c_2,\cA-\min(a_1,b_1,c_1),
\cB-\min(a_2,b_2,c_2)).$$
Let $\delta = 1$ if $k_0^{max} \geq k_0^{min}$ and $\cA,\cB\in\boZ$,
and let $\delta = 0$ otherwise. Denote the consecutive integers
$k_0^{min}$, $k_0^{min}+1$,..., $k_0^{max}$ by $k_0^{(i)}$, for
$1\leq i\leq M$. Then we have the tensor product multiplicity
$$M = Mult_{\lambda,\mu,\nu} = (k_0^{max} - k_0^{min} + 1)\delta$$
and the fusion rule coefficient on level $k$ is 
$$N_{\lambda,\mu,\nu}^{(k)} = \cases \max(i) &\text{ such that $k\geq k_0^{(i)}$
and $Mult_{\lambda,\mu,\nu} \neq 0$,} \cr \\
0 & \text{if $k < k_0^{(1)}$ or $Mult_{\lambda,\mu,\nu} = 0$}\endcases$$
so if $Mult_{\lambda,\mu,\nu} \neq 0$ and $k_0^{min}\leq k$ then
$N_{\lambda,\mu,\nu}^{(k)} = \min(k_0^{max},k) - k_0^{min} + 1$.
\epr

The fusion algebra for $A_2^{(1)}$ on level $k = 1$ is rather trivial, 
with three modules whose fusion rules are given exactly by the group
$\boZ_3$, which is the weight lattice modulo the root lattice of type $A_2$. 
The fusion algebra for $A_2^{(1)}$ on level $k=2$, and its covering by the 
group $\boZ_3^2$, are in the next two tables.
\vskip 10pt

\noindent{Table 8: Fusion Table for $A_2$ of level $k = 2$.}
\vskip 10pt

\begincellular{\centertable}
\row{}
\cell{[i]$\x$[j]}\cell{[0]}\cell{[1]}\cell{[2]}\cell{[3]}\cell{[4]}\cell{[5]}
\row{}
\cell{[0]}\cell{[0]}\cell{[1]}\cell{[2]}\cell{[3]}\cell{[4]}\cell{[5]}
\row{}
\cell{[1]}\cell{}\cell{[2]}\cell{[0]}\cell{[4]}\cell{[5]}\cell{[3]}
\row{}
\cell{[2]}\cell{}\cell{}\cell{[1]}\cell{[5]}\cell{[3]}\cell{[4]}
\row{}
\cell{[3]}\cell{}\cell{}\cell{}\cell{[0]+[3]}\cell{[1]+[4]}\cell{[2]+[5]}
\row{}
\cell{[4]}\cell{}\cell{}\cell{}\cell{}\cell{[2]+[5]}\cell{[0]+[3]}
\row{}
\cell{[5]}\cell{}\cell{}\cell{}\cell{}\cell{}\cell{[1]+[4]}
\endcellular

\newpage

\noindent{Table 9: Group $\boZ_3^2$ covering the Fusion Table 
for $A_2$ of level $k = 2$.}
\vskip 5pt

\begincellular{\centertable}
\row{}
\cell{[a]$\x$[b]}
\cell{(0,0)}\cell{(1,1)}\cell{(2,2)}
\cell{\vbox{\hbox{(1,2)} \hbox{(2,1)}} }
\cell{\vbox{\hbox{(2,0)} \hbox{(0,2)}} }
\cell{\vbox{\hbox{(1,0)} \hbox{(0,1)}} }
\row{}
\cell{(0,0)}
\cell{(0,0)}\cell{(1,1)}\cell{(2,2)}
\cell{\vbox{\hbox{(1,2)} \hbox{(2,1)}} }
\cell{\vbox{\hbox{(2,0)} \hbox{(0,2)}} }
\cell{\vbox{\hbox{(1,0)} \hbox{(0,1)}} }
\row{}
\cell{(1,1)}
\cell{}\cell{(2,2)}\cell{(0,0)}
\cell{\vbox{\hbox{(2,0)} \hbox{(0,2)}} }
\cell{\vbox{\hbox{(0,1)} \hbox{(1,0)}} }
\cell{\vbox{\hbox{(2,1)} \hbox{(1,2)}} }
\row{}
\cell{(2,2)}
\cell{}\cell{}\cell{(1,1)}
\cell{\vbox{\hbox{(0,1)} \hbox{(1,0)}} }
\cell{\vbox{\hbox{(1,2)} \hbox{(2,1)}} }
\cell{\vbox{\hbox{(0,2)} \hbox{(2,0)}} }
\row{}
\cell{\vbox{\hbox{(1,2)} \hbox{(2,1)}}}
\cell{}\cell{}\cell{}
\cell{\vbox{\hbox{(2,1),(0,0)} \hbox{(0,0),(1,2)}}}
\cell{\vbox{\hbox{(0,2),(1,1)} \hbox{(1,1),(2,0)}}}
\cell{\vbox{\hbox{(2,2),(0,1)} \hbox{(1,0),(2,2)}}}
\row{}
\cell{\vbox{\hbox{(2,0)} \hbox{(0,2)}}}
\cell{}\cell{}\cell{}\cell{}
\cell{\vbox{\hbox{(1,0),(2,2)} \hbox{(2,2),(0,1)}}}
\cell{\vbox{\hbox{(0,0),(1,2)} \hbox{(2,1),(0,0)}}}
\row{}
\cell{\vbox{\hbox{(1,0)} \hbox{(0,1)}}}
\cell{}\cell{}\cell{}\cell{}\cell{}
\cell{\vbox{\hbox{(2,0),(1,1)} \hbox{(1,1),(0,2)}}}
\endcellular
\vskip 5pt

In the fusion table for $A_2$ of level $k=3$, part of which is shown
below, we found a fusion coefficient greater than one, 
which our previous scheme could not handle. Our 
main problem has been to modify our method of covering fusion tables 
so as to account for such higher multiplicities. In doing so we found 
a surprising connection to a combinatorial problem which may be of 
some independent interest. In the next section we will set up the 
notation needed for this new point of view for $\hg$ of type 
$A_{N-1}^{(1)}$, and see how the combinatorics enters the picture.
\vskip 5pt

\noindent{Table 10: Partial fusion Table for $A_2$ of level $k = 3$.}
\vskip 5pt

\begincellular{\centertable}
\row{}
\cell{[i]$\x$[j]}\cell{[0]}\cell{[1]}\cell{[2]}\cell{[9]}
\row{}
\cell{[0]}\cell{[0]}\cell{[1]}\cell{[2]}\cell{[9]}
\row{}
\cell{[1]}\cell{}\cell{[2]}\cell{[0]}\cell{[9]}
\row{}
\cell{[2]}\cell{}\cell{}\cell{[1]}\cell{[9]}
\row{}
\cell{[9]}\cell{}\cell{}\cell{}\cell{[0]+[1]+[2]+2[9]}
\endcellular
\vskip 5pt

\noindent{\bf{3. A New Approach Incorporating Higher Multiplicities}}
\sk1

Let $\hg$ be the affine algebra of type $A_{N-1}^{(1)}$ built from 
$\bg = sl_N$. Let $\lambda_1,\cdots,\lambda_{N-1}$ denote the
fundamental weights of $\bg$, and let 
$$P^+ = \{\sum_{j=1}^{N-1} a_j\lambda_j\ |\ 0\leq a_j\in\boZ\}$$ 
denote the dominant integral weights of $\bg$. The irreducible 
modules $\hV_\Lambda$ for $\hg$ of level $k \geq 1$ are indexed by 
a dominant integral highest weight $\Lambda = kc + \lambda$ where 
$$\lambda = \sum_{j=1}^{N-1} a_j\lambda_j \in P^+$$ 
satisfies the ``level $k$ condition'' $\sum_{j=1}^{N-1} a_j \leq k$.
On level 1, those dominant integral weights are the fundamental
weights of $\hg$, 
$$\Lambda_0 = c, \Lambda_1 = c+\lambda_1,\cdots,\Lambda_{N-1} = 
c+\lambda_{N-1}.$$
Then the irreducible modules on level $k$ are in one-to-one 
correspondence with the set of $N$-tuples of nonnegative integers
$(i_0,i_1,\cdots,i_{N-1})$ whose sum is $k$. Such an $N$-tuple
corresponds to
$$\eqalign{\Lambda &= \sum_{j=0}^{N-1} i_j \Lambda_j \cr
&=\left(\sum_{j=0}^{N-1} i_j\right)c + i_1\lambda_1 + \cdots +
i_{N-1} \lambda_{N-1} \cr
&= kc + i_1\lambda_1 + \cdots + i_{N-1} \lambda_{N-1}.\cr}$$

Let $V_\lambda$ be the irreducible finite dimensional $g$-submodule
of $\hV_\Lambda$ generated by a highest weight vector. In the special
case when $\Lambda = k\Lambda_0 = kc$, that finite dimensional 
$\bg$-module is $V_0$, the one dimensional trivial $\bg$-module. 
Since $\bg$ is semisimple, any finite dimensional $\bg$-module is 
completely reducible. Therefore, we can write the tensor product 
of irreducible $\bg$-modules 
$$V_\lambda\otimes V_\mu = \sum_{\nu\in P^+} 
Mult_{\lambda,\mu}^\nu V_\nu$$
as the direct sum of irreducible $\bg$-modules, including multiplicities. 
This decomposition is independent of the level $k$ and is part of the 
basic representation theory of $sl_N$. The fusion algebra for $\hg$ on
level $k$ has a basis in one-to-one correspondence with the irreducible
$\hg$-modules on level $k$, and its fusion rules are obtained by a 
rather subtle truncation of the above summation. 

The fusion algebra for $A_{N-1}^{(1)}$ on level $k = 1$ is rather trivial, 
with a basis of $N$ vectors labelled by the highest weights of the $N$ 
irreducible fundamental modules, and with fusion rules given by the 
group $\boZ_N$, which is the weight lattice modulo the root lattice of 
$\bg$. This means that $N_{\lambda,\mu}^\nu = \delta_{\lambda_+\mu,\nu}$
where the addition takes place in the quotient group of the weight
lattice modulo the root lattice. 

Let $G$ be the group $\boZ_N^k$, and for each $N$-tuple of nonnegative
integers 
$$(i_0,i_1,\cdots, i_{N-1}) \hbox{ such that }
i_0+i_1+\cdots + i_{N-1} = k,$$ 
define the subset of $G$ 
$$P(i_0,i_1,\hdots,i_{N-1}) = \{x\in\boZ_N^k\ |\ j\hbox{ occurs exactly }
i_j\hbox{ times in }x,\ 0\leq j\leq N-1\}.$$
Then the cardinality of the subset $P(i_0,i_1,\hdots,i_{N-1})$ is
$${k\choose i_0,i_1,\hdots,i_{N-1}}$$
and $G$ is the disjoint union of these subsets. 
For example, with $N = 3$ and $k=2$, a bijection between these 
subsets and a basis of the $sl_3$ fusion algebra on level 2 is:
$$\eqalign{&[0] = \{(0,0)\} = P(2,0,0) \cr
&[1] = \{(1,1)\} = P(0,2,0) \cr
&[2] = \{(2,2)\} = P(0,0,2) \cr
&[3] = \{(1,2), (2,1)\} = P(0,1,1) \cr
&[4] = \{(2,0), (0,2)\} = P(1,0,1) \cr
&[5] = \{(1,0), (0,1)\} = P(1,1,0) \cr}$$

We can now start to describe our new method of covering 
the fusion algebra $F = F(g,k)$ for $\bg = sl_N$ on 
level $k$ by the group $G = \boZ_N^k$. First note that
the symmetric group $S_k$ acts on $G$ by permuting the $k$-tuples.
For $a\in G$, let $\cO_a = [a]$ denote the orbit of $a$ under this action
and let $\cO$ be the set of all such orbits. 
It is easy to see that these orbits are precisely the subsets
$P(i_0,i_1,\hdots,i_{N-1})$ defined above, and that each orbit
contains a unique representative in the ``standard'' form
$$(0^{i_0},1^{i_1},\cdots,(N-1)^{i_{N-1}})$$
where the exponent indicates the number of repetitions of the base. 
For $[a],[b],[c]\in \cO$ we wish to define a fusion coefficient
$N_{[a],[b]}^{[c]}$ in terms of the group $G$ so that 
the resulting product
$$[a] \cdot [b] = \sum_{[c]\in\cO} N_{[a],[b]}^{[c]} [c]$$
defines a fusion algebra isomorphic to the fusion algebra $F$. 
It is clear that the distinguished element $\Omega$ of the fusion algebra
should correspond to $[0]$ and that the conjugate of $[c]$ should be
$[-c]$. Then we can understand the fusion coefficients by studying
the completely symmetric coefficients 
$$N_{[a],[b],[c]} = N_{[a],[b]}^{[-c]}.$$
We may assume that 
$$\eqalign{
a &= (0^{\alpha_0},1^{\alpha_1},\cdots,(N-1)^{\alpha_{N-1}}) \cr
b &= (0^{\beta_0},1^{\beta_1},\cdots,(N-1)^{\beta_{N-1}}) \cr
c &= (0^{\gamma_0},1^{\gamma_1},\cdots,(N-1)^{\gamma_{N-1}}) \cr}$$
are in standard form. The exponents are nonnegative integers satisfying
$$\sum_{i=0}^{N-1} \alpha_i = \sum_{i=0}^{N-1} \beta_i =
\sum_{i=0}^{N-1} \gamma_i = k.$$

The group $S_k$ acts on
$$T([a],[b],[c]) = \{(x,y,z)\in [a]\x [b]\x [c] \ |\ x + y + z = 0 \}$$
which decomposes into a finite number of orbits under that action. 
We wish to determine that number of orbits, and show how it is related 
to the fusion coefficient $N_{[a],[b],[c]}$. We will show how to do this
for $N = 2$ and for $N = 3$ in the next two sections, and we hope that
future work will give the relationship for all ranks. 

\pr{Definition} Let $M([a],[b],[c])$ be the number of orbits of 
$T([a],[b],[c])$ under the action of $S_k$. \epr

Each orbit of $T([a],[b],[c])$ under $S_k$ has representatives whose third
component is in the standard form, and in fact, we may assume that 
$z = c$ is in that standard form. This makes it clear that the above counting
problem is equivalent to counting the number of orbits of 
$$T_c([a],[b]) = \{(x,y)\in [a]\x [b] \ |\ x + y + c = 0 \}$$
under the action of the stabilizer of $c$ in $S_k$,
$Stab_{S_k}(c) = \{\sigma\in S_k\ |\ \sigma(c)=c\}$. 
The structure of $Stab_{S_k}(c)$ is obviously the product of symmetric
groups $S_{\gamma_0}\x S_{\gamma_1}\x \cdots \x S_{\gamma_{N-1}}$ since
each subset of $\gamma_j$ occurences of $j$ can be freely permuted
without changing $c$. We now write   
$$x = (x_0,x_1,\dots,x_{N-1})\quad \hbox{ and }\quad
y = (y_0,y_1,\dots,y_{N-1})$$
where $x_i$ is the part of $x$ in the positions where $c$ has $i$'s and
where $y_i$ is the part of $y$ in the positions where $c$ has $i$'s.
By the action of the stabilizer of $c$, we can find a unique representative
of the orbit of that stabilizer so that each 
$$x_i = (0^{k_{i0}}, 1^{k_{i1}}, \dots, (N-1)^{k_{i,N-1}})$$
is in standard form. This defines the nonnegative integers $k_{ij}$ 
for $0\leq i,j\leq N-1$, which we use to form the matrix $K = [k_{ij}]$.
Then by the sum property, $x_i + y_i + c_i = (0^{\gamma_i})$, we must have
$$y_i = ( (N-i)^{k_{i0}}, (N-1-i)^{k_{i1}}, \dots, (1-i)^{k_{i,N-1}} ).$$
The $x$ and $y$ that have these forms are a system of distinct
representatives of the orbits, thus the number of such pairs $(x,y)$ is the
number $M([a],[b],[c])$. Below is a picture of the three vectors 
$x_i$, $y_i$ and $c_i$.
$$\matrix 
x_i\ = \ (\ &0^{k_{i0}}&,&  1^{k_{i1}}&,\ \dots\ ,&(N-1)^{k_{i,N-1}} &) \\
y_i\ = \ (\ &(N-i)^{k_{i0}}&,&  (N-1-i)^{k_{i1}}&,\ \dots\ ,&(1-i)^{k_{i,N-1}} &) \\
c_i\ = \ (\ &i^{k_{i0}}&,& i^{k_{i1}} &,\ \dots\ ,& i^{k_{i,N-1}} &) \endmatrix$$

The $3N$ equations which must be satisfied by the $N^2$ unknowns, 
$0\leq k_{ij} \in \boZ$, are
$$\eqalign{
&\sum_{j=0}^{N-1} k_{ij}  = \gamma_i,\qquad \hbox{for } 0\leq i\leq N-1 \cr
&\sum_{i=0}^{N-1} k_{ij}  = \alpha_j,\qquad \hbox{for } 0\leq j\leq N-1  \cr  
&\sum_{i=0}^{N-1} k_{i,l-i} = \beta_{N-l},\qquad \hbox{for } 1\leq l\leq N.\cr}$$
The first two equations are easy to understand from the definitions 
above as specifications of the sums in each row and column. 
But the third one, which specifies the sum of each upright diagonal 
(wrapped as if the matrix were toroidal), may need a little more 
consideration. When we write out the entries
$$\matrix 
y_0\ = &(N)^{k_{00}}&  (N-1)^{k_{01}}&\dots&(1)^{k_{0,N-1}} \\
y_1\ = &(N-1)^{k_{10}}&  (N-2)^{k_{11}}&\dots&(0)^{k_{1,N-1}} \\
y_2\ = &(N-2)^{k_{20}}&  (N-3)^{k_{21}}&\dots&(N-1)^{k_{2,N-1}} \\
\vdots \\
y_{N-1}\ = &(1)^{k_{N-1,0}}&  (0)^{k_{N-1,1}}&\dots&(2)^{k_{N-1,N-1}} \\
\endmatrix$$
we see that the base of the exponent $k_{ij}$ is $N-(i+j)$ modulo 
$N$, so the total number of occurences of base $N-l$ in $y$ is that
upright diagonal sum.

There are $N^2$ variables, and there are $3N$ equations, but the 
consistency of sums of different types means that 2 equations (of 
different types, say one column and one diagonal equation) 
are redundant. Therefore the dimension of the solution space 
(over $\bQ$) is $N^2 - (3N-2) = (N-1)(N-2)$, which is $0$
for $N=2$, $2$ for $N=3$, and $6$ for $N=4$.
\vskip 5pt

\noindent{\bf{4. The case of $N = 2$.}}
\sk1

For $N = 2$ with
$$a = (0^{\alpha_0},1^{\alpha_1}) \qquad
b = (0^{\beta_0},1^{\beta_1}) \qquad
c = (0^{\gamma_0},1^{\gamma_1})$$
and
$$\alpha_0 + \alpha_1 = \beta_0 + \beta_1 = \gamma_0 + \gamma_1 = k$$
we have the equations 
$$\matrix
k_{00} + k_{01} = \gamma_0 \qquad  &  k_{10} + k_{11} = \gamma_1 \\
k_{00} + k_{10} = \alpha_0 \qquad  &  k_{01} + k_{11} = \alpha_1 \\
k_{01} + k_{10} = \beta_1 \qquad  &  k_{00} + k_{11} = \beta_0 \endmatrix.$$
Three of these involving $k_{00}$ allow us to write
$$K = \bmatrix k_{00} & k_{01} \\ k_{10} & k_{11} \endbmatrix
= \bmatrix k_{00} & \gamma_0 - k_{00} \\ \alpha_0 - k_{00} & 
\beta_0 - k_{00} \endbmatrix,$$
so if a solution exists, then $k_{00}$ determines the matrix 
uniquely. Each of the other equations is equivalent to 
$$2k_{00} = \alpha_0 + \beta_0 + \gamma_0 - k 
= 2k - (\alpha_1 + \beta_1 + \gamma_1)$$
so
$$0\leq k_{00} = k - \hf(\alpha_1 + \beta_1 + \gamma_1) \in\boZ$$
imposes the requirement that
$$\hf(\alpha_1 + \beta_1 + \gamma_1) \in\boZ.$$
Since all entries of the matrix are nonnegative, we have
$0\leq k_{00}\leq min(\alpha_0,\beta_0,\gamma_0)$ so 
$$0\leq \alpha_0 + \beta_0 + \gamma_0 - k 
\leq 2\ min(\alpha_0,\beta_0,\gamma_0)$$
which is equivalent to
$$\max(2\alpha_1,2\beta_1,2\gamma_1)\leq 
\alpha_1 + \beta_1 + \gamma_1 \leq 2k.$$
Rewriting this as
$$\max(\alpha_1-\beta_1,\beta_1-\alpha_1,2\gamma_1-(\alpha_1+\beta_1))
\leq \gamma_1 \leq 2k -(\alpha_1+\beta_1)$$
and using $2\gamma_1-(\alpha_1+\beta_1) \leq \gamma_1$ iff 
$\gamma_1 \leq \alpha_1+\beta_1$, we get the condition
$$\max(\alpha_1-\beta_1,\beta_1-\alpha_1)
\leq \gamma_1 \leq \min(2k-(\alpha_1+\beta_1),\alpha_1+\beta_1)$$
which is equivalent to
$$| \alpha_1-\beta_1 | \leq \gamma_1 \leq \alpha_1+\beta_1 \qquad
\hbox{and}\qquad \alpha_1+\beta_1+\gamma_1 \leq 2k.$$

\pr{Theorem 3} With notation as above, for $1\leq k\in \boZ$,
we have $M([a],[b],[c]) = 1$ if 
$$| \alpha_1-\beta_1 | \leq \gamma_1 \leq \alpha_1+\beta_1,\quad
\alpha_1+\beta_1+\gamma_1 \leq 2k\quad\hbox{and}\quad
\hf(\alpha_1 + \beta_1 + \gamma_1) \in\boZ,$$
$M([a],[b],[c]) = 0$ otherwise. \epr

\pr{Corollary} When $N = 2$, with notation as above, and with
$a_1 = \alpha_1/2$, $b_1 = \beta_1/2$ and $c_1 = \gamma_1/2$, we have 
$M([a],[b],[c]) = N_{[a_1],[b_1],[c_1]}$.
\epr
\demo{Proof} The element $a = (0^{\alpha_0},1^{\alpha_1}) \in\boZ_2^k$ 
corresponds to a level $k$ highest weight $\Lambda = \alpha_0\Lambda_0 +
\alpha_1\Lambda_1 = kc + \alpha_1\lambda_1$, for $0\leq\alpha_1\leq k$
integral. The associated irreducible highest weight $sl_2$-module is
then $V_{\alpha_1/2}$ if we use the spin $a_1 = \alpha_1/2$ to label
the module. The labels for the $sl_2$-modules associated with $b$ and $c$
are $b_1 = \beta_1/2$ and $c_1 = \gamma_1/2$. Since $c = -c\in\boZ_2^k$,
$N_{[a_1],[b_1],[c_1]} = N_{[a_1],[b_1]}^{[c_1]}$, and the conditions
in the theorem above translate exactly into the conditions of the fusion
rules for $sl_2$. \enddemo
\vskip 5pt

\noindent{\bf{5. The case of $N = 3$.}}
\sk1

For $N = 3$, after using two column sums, two row sums and one diagonal
sum to rewrite the entries in the last row and column of $K$, we have
$$K = \bmatrix k_{00} &  k_{01} &  \gamma_0-k_{00}-k_{01}  \\
k_{10} &  k_{11}& \gamma_1-k_{10}-k_{11}  \\
\alpha_0-k_{00}-k_{10}  &  \alpha_1-k_{01}-k_{11} & \beta_2-k_{01}-k_{10}
\endbmatrix.$$
This leaves 4 more equations:
$$\eqalign{
[\alpha_0-k_{00}-k_{10}] + [\alpha_1-k_{01}-k_{11}] + [\beta_2-k_{01}-k_{10}] &= \gamma_2,  \cr
[\gamma_0-k_{00}-k_{01}] + [\gamma_1-k_{10}-k_{11}] + [\beta_2-k_{01}-k_{10}] &= \alpha_2,  \cr
k_{00} + [\gamma_1-k_{10}-k_{11}] + [\alpha_1-k_{01}-k_{11}] &= \beta_0,  \cr
[\gamma_0-k_{00}-k_{01}] + k_{11} + [\alpha_0-k_{00}-k_{10}] &= \beta_1.\cr}$$
Using the relation 
$$k = \alpha_0 + \alpha_1 + \alpha_2 = \beta_0 + \beta_1 + \beta_2 =
\gamma_0 + \gamma_1 + \gamma_2$$
we find that two of the four equations are redundant, and the system is
equivalent to 
$$\eqalign{
k_{00} + 2k_{01} + 2k_{10} + k_{11} &= k - \alpha_2 + \beta_2 - \gamma_2, \cr
3(k_{00} - k_{11}) &= (\alpha_0+\beta_0+\gamma_0)-(\alpha_1+\beta_1+\gamma_1).\cr}$$
Letting 
$$A = \td[(\alpha_0+\beta_0+\gamma_0)-(\alpha_1+\beta_1+\gamma_1)]$$
the solutions to the linear system above are
$$k_{00} = k_{11} + A \qquad k_{01} = -k_{10}-k_{11}+A+\alpha_1-\beta_0+\gamma_1$$
where $k_{10}$ and $k_{11}$ are the ``free variables''. In fact, because
the entries of the matrix $K$ are nonnegative integers, we have $A\in\boZ$ and
the following nine inequalities:
$$0\leq k_{00} \quad 0\leq k_{01} \quad 0\leq k_{10} \quad 0\leq k_{11} 
\quad k_{01}+k_{10}\leq \beta_2$$
$$k_{00}+k_{10}\leq \alpha_0 \quad k_{01}+k_{11}\leq \alpha_1 \quad 
k_{00}+k_{01}\leq \gamma_0 \quad k_{10}+k_{11}\leq \gamma_1.$$
These translate into the following inequalities on the free variables 
$k_{10}$ and $k_{11}$:
$$-A \leq k_{11} \qquad k_{10}+k_{11}\leq A+\alpha_1-\beta_0+\gamma_1 \qquad
0\leq k_{10} \qquad 0\leq k_{11}$$
$$k_{10}+k_{11}\leq \alpha_0-A \qquad A-\beta_0+\gamma_1 \leq k_{10} 
\qquad 2A+\alpha_1-\beta_0+\gamma_1-\gamma_0\leq k_{10}$$
$$k_{10}+k_{11}\leq \gamma_1 \qquad A+\alpha_1-\beta_0-\beta_2+\gamma_1\leq k_{11}.$$
Combining these, we can say that
$$C = \max(0,A-\beta_0+\gamma_1,2A+\alpha_1-\beta_0+\gamma_1-\gamma_0) \leq k_{10}$$
$$D = \max(-A,0,A+\alpha_1-\beta_0-\beta_2+\gamma_1) \leq k_{11}$$
$$k_{10}+k_{11}\leq \min(A+\alpha_1-\beta_0+\gamma_1,\alpha_0-A,\gamma_1) = E.$$
It is clear that when $E-C \geq D$ and $A\in\boZ$,
these constraints define a ``45-45-90'' triangle in the first 
quadrant of the $k_{10} k_{11}$-plane whose vertices have integral coordinates.
The number $M([a],[b],[c])$ is exactly the number of integral points within 
that triangle, including its boundary. The number of points on one leg of 
the triangle is clearly $E-D-C+1$ so that $M([a],[b],[c])$ is a triangular number.

\pr{Theorem 4} With notation as above, for $1\leq k\in \boZ$,
we have
$$M([a],[b],[c]) = \frac{(E-D-C+1)(E-D-C+2)}{2} = {E-D-C+2 \choose 2}.$$ 
\epr

We will now show that $E-D-C+1$ is equal to the fusion rule coefficient 
$N_{[a],[b],[c]}$ by using the theorem of \cite{BMW}. 
The level $k$ highest weights, $\lambda$, $\mu$ and $\nu$, referred to in the 
\cite{BMW} theorem, correspond to $a$, $b$ and $c$ in $\boZ_3^k$ so that
$$\lambda = \alpha_1 \lambda_1 + \alpha_2 \lambda_2\ , \qquad
\mu = \beta_1 \lambda_1 + \beta_2 \lambda_2 \ ,\qquad
\nu = \gamma_1 \lambda_1 + \gamma_2 \lambda_2.$$
We have
$$\cA = \td\left[{2(\alpha_1+\beta_1+\gamma_1) + (\alpha_2+\beta_2+\gamma_2)}\right],$$
$$\cB = \td\left[{(\alpha_1+\beta_1+\gamma_1) + 2(\alpha_2+\beta_2+\gamma_2)}\right],$$
so that
$$\cA = k - A \ ,\qquad 2\cA-\cB=\alpha_1+\beta_1+\gamma_1, \qquad 
2\cB-\cA=\alpha_2+\beta_2+\gamma_2,$$
and
$$\cA+\cB = (\alpha_1+\beta_1+\gamma_1) + (\alpha_2+\beta_2+\gamma_2).$$
The variables used in the \cite{BMW} theorem are
$$\eqalign{
k_0^{max} &= min(\cA,\cB) \cr
k_0^{min} &= max(\alpha_1+\alpha_2,\beta_1+\beta_2,\gamma_1+\gamma_2,
\cA-min(\alpha_1,\beta_1,\gamma_1),\cB - min(\alpha_2,\beta_2,\gamma_2)) \cr
&= max(k-min(\alpha_0,\beta_0,\gamma_0),\cA-min(\alpha_1,\beta_1,\gamma_1),
\cB - min(\alpha_2,\beta_2,\gamma_2)).\cr}$$
The \cite{BMW} theorem says that the fusion coefficient 
$$N_{\lambda,\mu,\nu}^{(k)} = \min(k_0^{max},k) - k_0^{min} + 1$$
if the conditions
$$k_0^{min} \leq k_0^{max}\ ,\quad k_0^{min} \leq k \quad\hbox{and}
\quad \cA,\cB\in\boZ$$
are satisfied, and the coefficient is zero otherwise. 
Using the formulas above, it is straightforward to check that
$$C = \max(0,\beta_1+\beta_2+\gamma_1-\cA,\cB-\cA-\alpha_2+\gamma_1),$$
$$D = \cA - \min(\cA,\cB,k) = \cA - \min(k_0^{max},k),$$
$$E = \cA - \max(\alpha_1+\alpha_2,\cA-\gamma_1,\cB-\beta_2).$$
Two identities useful for checking this are
$$2A+\alpha_1-\beta_0+\gamma_1-\gamma_0 = \cB - \cA - \alpha_2 + \gamma_1,$$
$$A+\alpha_1-\beta_0-\beta_2+\gamma_1 = \cB - (\alpha_2+\beta_2+\gamma_2).$$
{}From these expressions for $C$, $D$ and $E$ we can now compute
$$E-D-C+1 =\min(k_0^{max},k)+1-$$
$$\{\max(0,\beta_1+\beta_2-\cA+\gamma_1,\cB-\cA-\alpha_2+\gamma_1)+
\max(\alpha_1+\alpha_2,\cA-\gamma_1,\cB-\beta_2)\}$$
$$=\min(k_0^{max},k)+1-$$
$$\max(\alpha_1+\alpha_2, \cA-\gamma_1, \cB-\beta_2,\cB-\gamma_2,
\beta_1+\beta_2,\cA-\alpha_1,\cA-\beta_1, \cB-\alpha_2, \gamma_1+\gamma_2)$$
$$=\min(k_0^{max},k)+1-k_0^{min}.$$

We also have 
$$E-C=\cA-k_0^{min}$$
so the condition $E-C\geq D$ is equivalent to the condition
$$\min(k_0^{max},k) \geq k_0^{min},$$
that is,
$$k_0^{max} \geq k_0^{min}\hbox{  and  } k\geq k_0^{min}.$$
It is clear that $A \in \boZ$ iff $\cA \in \boZ$ iff $\cB \in \boZ$
since $\cA = k-A$ and $\cA + \cB \in\boZ$.
We have shown that the conditions for the fusion coefficient to be
nonzero are equivalent to the conditions for $M([a],[b],[c]) \neq 0$.
We have now proved our main result.

\pr{Theorem 5} For $N = 3$, the number of orbits of $T([a],[b],[c])$ under the 
action of $S_k$ is related to the fusion coefficient $N_{[a],[b],[c]}$ by
$$M([a],[b],[c]) = {N_{[a],[b],[c]} + 1 \choose 2}.$$
\epr
 
Another form of the conditions and number of orbits can be given
as follows. For $i,j,l,r\in\{0,1,2\}$ taken modulo $3$, let
$$s_{ijl} = \alpha_i + \beta_j + \gamma_l,$$
and
$$D^r_{ij} = s_{i,i+r,i+2r} - s_{j,j+r,j+2r}.$$
Then, using the fact that for $(i,j,l)$ any permutation of $(0,1,2)$, 
$$s_{i,i+r,i+2r} + s_{l,l+r,l+2r} = 3k - s_{j,j+r,j+2r},$$
we get that
$$D^r_{ij} - D^r_{jl} = 3(k - s_{j,j+r,j+2r})$$
so 
$$D^r_{ij} \equiv D^r_{jl} \pmod 3.$$
In fact, we have that
$D^0_{01} \equiv D^0_{12} \equiv D^0_{20} \equiv D^1_{01} \equiv D^1_{12} 
\equiv D^1_{20} \equiv D^2_{01} \equiv D^2_{12} \equiv D^2_{20}\pmod 3$.
This makes the necessary condition written below ($3|D$) 
for existence of a solution, independent of which 
$D^r_{ij}$ one checks for divisibility by $3$.

\pr{Theorem 6} [Zaslavsky] Let $D$ be any $D^r_{i,i+1}$ defined above and define $m$ by
$$3m = \min(s_{000}, s_{111}, s_{222}) + \min(s_{012}, s_{120}, s_{201}) 
+ \min(s_{021}, s_{102}, s_{210}) - 2k.$$
If $3\mid D$ and $m\geq 0$ then we have 
$m\in\boZ$ and $M([a],[b],[c]) = \binom{m+2}{2}$. 
Otherwise, $M([a],[b],[c]) = 0$. \epr
\sk1

\enddocument
\bye